\documentclass{ifacconf}

\usepackage{graphicx}      
\usepackage{natbib}        

\usepackage{amsmath,amssymb,amsfonts}
\usepackage{algorithmic}
\usepackage{graphicx}
\usepackage{textcomp}
\def\BibTeX{{\rm B\kern-.05em{\sc i\kern-.025em b}\kern-.08em
    T\kern-.1667em\lower.7ex\hbox{E}\kern-.125emX}}
\usepackage{url}
\usepackage[english]{babel}
\usepackage[utf8]{inputenc}
\usepackage{amsmath}
\usepackage{amsfonts}
\usepackage{amssymb}
\usepackage{mathtools}
\usepackage{amsmath}
\usepackage{color}
\usepackage{listings}
\usepackage{setspace}
\usepackage{nicefrac}
\usepackage{csquotes}
\usepackage{calc}
\usepackage{lipsum}
\usepackage{url}
\usepackage{longtable}
\usepackage{stackrel}
\usepackage{ mathrsfs }
\usepackage{matlab-prettifier}
\usepackage{soul}
\usepackage{enumitem}
\usepackage{tabularx}
\usepackage{bbm}
\usepackage{extarrows}
\usepackage[normalem]{ulem}
\usepackage[official]{eurosym}
\usepackage{booktabs}
\usepackage{flushend}
\newcommand{\R}{\mathbb{R}}

\newcommand{\Z}{\mathbb{Z}}

\newcommand{\B}{\mathbb{B}}

\newcommand{\E}{\mathbb{E}}
\renewcommand{\P}{\mathbb{P}}

\newcommand{\td}{\text{\normalfont d}}

\newcommand{\eqb}{\begin{flalign}}
\newcommand{\eqe}{\end{flalign}}
\newcommand{\Ra}{\Rightarrow}

\newcommand{\mc}{\mathcal}

\newcommand{\barb}[1]{\bar{\mathbf #1}}
\newcommand{\tmc}[1]{\widetilde{\mc #1}}

\newcommand{\cl}{\text{\normalfont cl}}
\newcommand{\omitcustom}[1]{{}}

\newcommand{\ediamond}[0]{\mbox{}\hfill$\Diamond$}

\newcommand{\TIMETEMP}{T}

\newlength{\leftstackrelawd}
\newlength{\leftstackrelbwd}
\def\leftstackrel#1#2{\settowidth{\leftstackrelawd}%
{${{}^{#1}}$}\settowidth{\leftstackrelbwd}{$#2$}%
\addtolength{\leftstackrelawd}{-\leftstackrelbwd}%
\leavevmode\ifthenelse{\lengthtest{\leftstackrelawd>0pt}}%
{\kern-.5\leftstackrelawd}{}\mathrel{\mathop{#2}\limits^{#1}}}

\makeatletter
\let\old@ssect\@ssect 
\makeatother
\usepackage{natbib}
\usepackage[hidelinks]{hyperref}
\makeatletter
\def\@ssect#1#2#3#4#5#6{%
  \NR@gettitle{#6}
  \old@ssect{#1}{#2}{#3}{#4}{#5}{#6}
}
\makeatother

\begin{document}
\begin{frontmatter}

\title{Robust Recurrence of Discrete-Time Infinite-Horizon Stochastic Optimal Control with Discounted Cost\thanksref{footnoteinfo}}

\thanks[footnoteinfo]{This work was supported by the ARC under the Discovery Project DP210102600, by the France-Australia collaboration project IRP-ARS CNRS,  and by the ANR grant OLYMPIA ANR-23-CE48-0006.}



\author[First,Second]{Robert H. Moldenhauer}
\author[First]{Dragan Nešić}
\author[First]{Mathieu Granzotto}
\author[Second]{Romain Postoyan}
\author[Third]{Andrew R. Teel}
\address[First]{Department of Electrical and Electronic Engineering, University of Melbourne, Parkville, VIC 3010, Australia (e-mail: moldenhauer.r@student.unimelb.edu.au, \{dnesic,mathieu.granzotto\}@unimelb.edu.au)}
\address[Second]{Université de Lorraine, CNRS, CRAN, F-54000 Nancy, France (e-mail: \{robert.moldenhauer, romain.postoyan\}@univ-lorraine.fr)}
\address[Third]{Department of Electrical and Computer Engineering, University of California, Santa Barbara, CA 93106 USA (e-mail: teel@ece.ucsb.edu)}

\begin{abstract}~                
We analyze the stability of general nonlinear discrete-time stochastic systems controlled by optimal inputs that minimize an infinite-horizon discounted cost.
Under a novel stochastic formulation of cost-controllability and detectability assumptions inspired by the related literature on deterministic systems, we prove that uniform semi-global practical recurrence holds for the closed-loop system, where the adjustable parameter is the discount factor.
Under additional continuity assumptions, we further prove that this property is robust.
\end{abstract}


\end{frontmatter}

\section{Introduction}
Optimal control seeks to minimize a cost incurred during the system operation \citep{kirk_2004_optimal}, which is a naturally compelling way to ensure desirable performance properties for the controlled system.
Optimal control is able to address decision-making problems not only in control engineering, but also in artificial intelligence, operations research, economics, medicine, and so on.
We talk of stochastic optimal control when the dynamical system includes a random input and expected cost is minimized.
This occurs in a range of applications involving substantial uncertainty, such as electric vehicles \citep{moura_stochastic_2011}, wireless networked control systems \citep{wang2016stochastic}, biological systems \citep{ulas_acikgoz_blood_2010, todorov_stochastic_2005}, and finance \citep{soner2004stochastic}.
Moreover, reinforcement learning is often formulated in terms of such stochastic optimal control problems, usually referred to as Markov decision processes in that context \citep{sutton_reinforcement_2018}.

In this paper, we focus on stochastic optimal control in discrete time over an infinite horizon with discounted costs.
The optimal control problem is solved over stationary policies, and discounting ensures that policies with finite costs exist for a large class of stochastic systems.
Our objective is to derive conditions under which the closed-loop feedback system resulting from near-optimal policies is stable.
Allowing near-optimal, not just optimal policies is of practical importance, as one may not be able to find exact optimal policies algorithmically, and approximate dynamic programming is employed instead \cite{bertsekas_dynamic_2012}.
The recent work by \cite{granzotto_stability_2024} achieved stability guarantees for the stochastic linear discounted quadratic problem.
We are not aware of any equivalent results for general nonlinear systems and costs, which are therefore the subject of this paper.


\omitcustom{While stability of the closed loop is of vital importance in many control applications, such is robustness to model perturbations.
In particular, optimality often relies on precise knowledge of the system and analytical solutions, while we may only have access to approximate dynamics and approximate or near-optimal inputs.
Under certain conditions, stability implies existence of a continuous Lyapunov function \cite{kellett_discrete-time_2004}, which in turn certifies robust stability \cite{kellett_smooth_2004,kellett_robustness_2005}.
This is also the case for stochastic systems, see \cite{subbaraman_converse_2013, teel_converse_2014}.
With this in mind, a co-design approach that leverages both optimality and stability to build continuous Lyapunov functions can yield powerful properties of the closed-loop.
This has been successfully employed in the deterministic setting for nonlinear systems and costs, see \cite{grimm_model_2005,postoyan_stability_2017,granzotto_finite-horizon_2021,granzotto_robust_2024}.
There, it was shown that (robust) stability is not only achieved by optimal policies, but also by near-optimal inputs generated by approximate dynamic programming procedures.
We are not aware of any equivalent results for stochastic optimal control problems at the exception of the recent work in \cite{granzotto_stability_2024} dedicated to the stochastic LQR problem.
}

We found that the notion of \emph{semi-global practical recurrence}, which was first introduced by \cite{grammatico_discrete-time_2013}, best describes the general behaviour of the considered closed loop.
Roughly speaking, a set is said to be recurrent if almost all solutions revisit the set infinitely often, see \citep{subbaraman_converse_2013}, which is inspired by \cite{meyn_markov_1993}.
Semi-global practical recurrence is a weaker property, that \cite{grammatico_discrete-time_2013} motivated through the following analogy:
\emph{``A walking robot in a stochastically rough terrain is guaranteed
to eventually fall down (with probability one as $t\to\infty$).
However, one can still meaningfully quantify stochastic stability in terms of expected time to failure.
As long as the robot does not fall down, it keeps returning back to certain walking postures, that we would hence call recurrent''}
\cite[Section 7]{grammatico_discrete-time_2013}.
We consider this property with respect to the discount factor $\gamma\in(0,1)$.
It then means that, if $\gamma$ is sufficiently close to 1, failure occurs late in expectation and, until that happens, states return near the recurrent set with high probability.
For near-optimal instead of optimal policies, semi-global practical recurrence is preserved, albeit to a generally larger recurrent set depending on the cost gap.

Under additional conditions, we show robustness to model perturbations.
This is important, as the optimal policies are typically designed on an ideal plant model, which may differ from the real-world plant.
The semi-global practical recurrence is then with respect to the discount factor and a bound on the magnitude of disturbances.

The main tools used to obtain our results are Bellman equation and Lyapunov analysis.
This follows in the footsteps of existing results on deterministic systems by \cite{grimm_model_2005,postoyan_stability_2017,granzotto_finite-horizon_2021,granzotto_stability_2024,granzotto_robust_2024}. 
There, under certain cost-controllability and detectability assumptions, a Lyapunov function is constructed by combining the optimal value function and Lyapunov like function characterizing detectability, such that its decrease is guaranteed by combining the Bellman equation with detectability.
In this paper, we weaken the cost-controllability and detectability assumptions made by \cite{postoyan_stability_2017} in order to address a large class of stochastic systems.
We demonstrate that the new assumptions generalize the stochastic linear quadratic problem as a special case.
Under these assumptions as well as some mild regularity and measurability conditions, we prove the desired semi-global practical recurrence.
Robustness is then achieved under additional continuity conditions.

\omitcustom{The main tools we use to approach the intersection of optimality and stability are dynamic programming from the optimality side, and Lyapunov analysis for the stability side.
Under appropriate cost controllability and detectability assumptions, we construct a Lyapunov function by combining the optimal value function and Lyapunov like function characterizing detectability, where its decrease is guaranteed by combining the Bellman equation with detectability.
This technique was successfully employed on deterministic systems by \cite{grimm_model_2005} for model predictive control, by \cite{postoyan_stability_2017} for infinite-horizon discounted and by \cite{granzotto_finite-horizon_2021} for (in)finite-horizon (un)discounted, and linear stochastic setting in \cite{granzotto_stability_2024}.
We adapt this procedure to the stochastic nonlinear case with discounted infinite-horizon costs.
We weaken controllability and detectability assumptions to hold for a large class of stochastic systems, and demonstrate they generalize the linear quadratic problem as a special case.
Under these assumptions and some mild regularity and measurability conditions, we establish semi-global practical recurrence with respect to the discount factor.
Finally, if the system, stage cost and value function are continuous, we prove that this property is robust to small perturbations.}

The remainder of this paper is organized as follows. Section
2 contains notation and basic definitions. In Section 3,
we introduce the system and optimal control problem
and define strictly causal random closed-loop solutions.
We start Section 4 by introducing the stochastic versions
of cost-controllability and detectability assumptions and
illustrate them on an example, before stating the main
result on semi-global practical recurrence, whose proof can be found in the appendix.
Section 5 states the robust version under additional continuity assumptions.
Section 6 concludes the paper.

\omitcustom{
Optimal control selects control inputs so as to minimize a cost incurred during the system operation \cite{kirk_2004_optimal}.
It is able to address decision-making problems not only in control engineering, but also in artificial intelligence, operations research, economics, medicine, etc.
In this paper, we focus on stochastic optimal control in discrete time over an infinite horizon, with general costs and for general nonlinear system dynamics involving random inputs.
Reinforcement learning is often formulated in terms of such stochastic optimal control problems, usually known as Markov decision processes in that context \cite{sutton_reinforcement_2018}.
Stochastic optimal control is being researched for a wide range of applications involving substantial uncertainty, such as electric vehicles \cite{moura_stochastic_2011}, wireless networked control systems \cite{wang2016stochastic}, biological systems \cite{ulas_acikgoz_blood_2010, todorov_stochastic_2005} and finance \cite{soner2004stochastic}.
While stability in itself is crucial in most control applications, a co-design that targets optimality and stability together proved in the deterministic case to yield benefits such as robustness and near-optimality of numerical approximate dynamic programming methods \cite{grimm_model_2005,postoyan_stability_2017,granzotto_finite-horizon_2021}.
In this paper we prove robust semi-global practical recurrence for discounted stochastic nonlinear optimal control problems, intended to be a first step in unlocking the potential of such a co-design to stochastic systems.

The main tools used to approach the intersection of optimality and stability are dynamic programming from the optimality side, and Lyapunov analysis for stability.
Under appropriate controllability and detectability assumptions a Lyapunov function based on the optimal value function can be constructed, with its decrease guaranteed by the Bellman equation.
This technique was successfully employed on deterministic systems by \cite{grimm_model_2005} for finite-horizon undiscounted problems (i.e., model predictive control), by \cite{postoyan_stability_2017} for infinite-horizon discounted and by \cite{granzotto_finite-horizon_2021} for finite-horizon discounted.
The obtained form of stability is semi-global and practical in terms of the prediction horizon $N$, discount factor $\gamma$ or both, respectively.
\textit{Semi-global} means, that the domain of attraction gets larger with higher $N$ (or $\gamma$, or both), including every point eventually as $N\to\infty$ (or $\gamma\to1$, or both).
\textit{Practical} means, that the asymptotically stable set gets smaller with higher $N$ and $\gamma$, excluding every point outside a limit set eventually.

In this paper, we adapt this procedure to the stochastic case.
We weaken the controllability and detectability assumptions to hold for a large class of stochastic systems, and demonstrate them on the linear quadratic problem with additive noise.
Under these assumptions we prove semi-global practical recurrence with respect to the discount factor $\gamma$ as it approaches 1.
This novel notion of stability was first introduced by \cite{grammatico_discrete-time_2013}, who explained it through the following analogy: "A walking robot in a stochastically rough terrain is guaranteed to eventually fall down (with probability one as $t\to\infty$).
However, one can still meaningfully quantify stochastic stability in terms of expected time to failure. As long as the robot does not fall down, it keeps returning back to certain walking postures, that we would hence call recurrent" \cite[Section 7]{grammatico_discrete-time_2013}.
Semi-global practical recurrence is the property that, with $\gamma$ close to 1, failure occurs late in expectation and, until then, the state returns near the recurrent set with high probability.
Finally, if the system, stage cost and value function are continuous, we prove that this property is robust to small perturbations.

The remainder of this paper is organized as follows.
Section \ref{s:Notation} contains notation and basic definitions.
In Section \ref{s:setup}, we introduce the system and optimal control problem and define strictly causal random closed-loop solutions.
We start Section \ref{s:Recurrence} by introducing stochastic versions of cost-controllability and detectability assumptions and demonstrate them on an example, before stating the main result on semi-global practical recurrence.
Section \ref{s:Robust_Recurrence} states the robust version of this under additional continuity assumptions.
Section \ref{s:Conclusion} concludes the paper.
}

\section{Notation and basic definitions} \label{s:Notation}
    \omitcustom{\item} The symbol $\R$ denotes the set of real numbers, $\R_{\geq0}$ the set of non-negative real numbers and $\R_{>0}$ the set of positive real numbers.
    The notation $\Z_{\geq0}$ ($\Z_{>0}$, resp.) denotes the nonnegative (positive, resp.) integers.
    The Euclidean norm is denoted by $|\cdot|$.
    The spectral norm of a matrix is denoted by $||\cdot||_2$.
    \omitcustom{\item} Given sets $X,Y$, a function $f:X\to Y$ and $B\subseteq Y$, the pre-image of $B$ under $f$ is denoted by $f^{-1}(B)$.
    The notation $\B^n:=\{x\in\R^n~|~|x|<1\}$ denotes the $n$-dimensional open unit ball of $\R^n$ centered around the origin, where we just write $\B$ when $n$ is clear from context.
    \omitcustom{\item} For any set $\mc S\subseteq\R^n$, the notation $\cl(\mc S)$ denotes the closure of $\mc S$. 
    \omitcustom{\item} A function $\alpha:\R_{\geq0}\to\R_{\geq0}$ is said to \textit{belong to class} $\mc K_\infty$ ($\alpha\in\mc K_\infty$) if it is continuous, strictly increasing, zero at zero and unbounded.
    A function $f:X\to Y$, where $X\subseteq\R^n$ and $Y\subseteq\R^m$, is said to be \textit{locally bounded} if images of compact sets are bounded.
    A function $f:X\to\R_{\geq0}$, where $X\subseteq\R^n$, is said to be \textit{radially unbounded} if pre-images of bounded sets are bounded; note that if $X$ is bounded, then $f$ is always radially unbounded.
    A \textit{set-valued mapping} $\Gamma$ from a set $X$ to a set $Y$ is a mapping from $X$ to subsets of $Y$, and is notated as $\Gamma:X\rightrightarrows Y$.
    
    \omitcustom{\item} A \textit{measurable space} is a pair $(\Omega,\mc F)$ of a set $\Omega$ and a $\sigma$-algebra $\mc F$ over $\Omega$.
    Over subsets of $\R^n$, it is implied that the Borel $\sigma$-algebra is used as $\sigma$-algebra.
    A set $F\subseteq\Omega$ is said to be \textit{$\mc F$-measurable} (or just \textit{measurable}, if $\mc F$ is clear from context) if $F\in\mc F$.
    A \textit{probability space} $(\Omega,\mc F,\P)$ consists of a measurable space $(\Omega,\mc F)$ and a probability measure $\P:\mc F\to[0,1]$.
    \omitcustom{\item} Given measurable spaces $(\Omega,\mc F)$ and $(\mc X,\mc G)$, a function $f:\Omega\to\mc X$ is said to be \textit{measurable}, if $f^{-1}(B)\in\mc F$ for any $B\in\mc G$.
    If, in addition, $(\Omega,\mc F,\P)$ is a probability space, then such a measurable function is referred to as \textit{random variable}, and we denote it by a boldface lowercase letter.
    \omitcustom{\item} The expected value of a random variable $\mathbf x$ is denoted by $\E[\mathbf x]=\int \mathbf x\td\P=\int \mathbf x(\omega)\td\P(\omega)$.
    For $F\in\mc F$ (such a set is referred to as an \textit{event}) with $\P(F)>0$, the \textit{conditional expectation} of $\mathbf x$ is denoted by $\E[\mathbf x|F] = \P(F)^{-1}\int_F\mathbf x\td\P$.
    We use round brackets, i.e., $\P(F)$, for the probability of an event $F\in\mc F$, and square brackets, i.e., $\P[\cdot]$, on logical predicates involving random variables, as short-hand notation for the set of all $\omega\in\Omega$ that satisfy the logical expression.


    Whenever we apply a function or arithmetic operation on a set or a random variable, we understand it to be the set or random variable that results when the function or operation is applied element-wise or path-wise, respectively.
    For example, if $f$ is a function and $X$ is a set, $f(X):=\{f(x)~|~x\in X\}$. If $f$ is a function and $\mathbf x:\Omega\to\mc X$ is a random variable, then $f(\mathbf x)$ is the random variable with $f(\mathbf x)(\omega) = f(\mathbf x(\omega))$ for all $\omega\in\Omega$.
\section{Problem Statement}\label{s:setup}
In this section, we introduce the system and costs, followed by defining a set-valued near-optimal feedback policy and the strictly causal random closed-loop solutions it generates.
\subsection{Class of systems and costs}\label{s:3.1:Classofsystemsandcosts}
Consider a function $f:\mc X\times\mc U\times\mc V\to\mc X$, where $\mc X\subseteq\R^n$ and $\mc U\subseteq\R^m$ are closed, $\mc V\subseteq\R^q$ is measurable and $n,m,q\in\Z_{>0}$, and the discrete-time stochastic control system
\begin{align}
    x^+ = f(x,u,v)\label{eq:nonlinear_system}
\end{align}
with state $x\in\mc X$, control input $u\in\mc U$ and random input $v\in\mc V$.
Let a probability space $(\Omega,\mc F,\P)$, a random variable $\mathbf v:\Omega\to\mc V$ and an independently identically distributed (i.i.d.) sequence of random variables $(\mathbf v_0, \mathbf v_1, \dots)$
with the same distribution as $\mathbf v$ be given, where $\mathbf v_k$ acts as random input at time $k$.


In optimal control of stochastic systems, the optimization problem is often solved over feedback policies, rather than sequences of inputs (which is customary for deterministic optimal control, e.g., \cite{grimm_model_2005, postoyan_stability_2017, granzotto_finite-horizon_2021}).
Policies allow the controller to respond to random inputs through their effect on the current state, see, e.g., \cite{Hern_ndez_Lerma_1996, bertsekas_dynamic_2012}.
Let the set of admissible (feedback) policies be
\begin{align}
    \mc H_a:=\{h:\mc X\to\mc U~|~h\text{ is measurable}\}.
\end{align}
We denote the random closed-loop solution to \eqref{eq:nonlinear_system} at the $k^\text{th}$-step starting at state $x\in\mc X$ and under feedback policy $h\in\mc H_a$ by $\boldsymbol\varphi(k,x,h):\Omega\to\mc X$, that is,
\begin{align}
\begin{split}
    \boldsymbol\varphi(k+1,x,h) &:= f(\boldsymbol\varphi(k,x,h), h(\boldsymbol\varphi(k,x,h)), \mathbf v_k),\label{phi_definition}\\
    \boldsymbol\varphi(0,x,h) &:= x.
\end{split}
\end{align}
Note that, as $\mathbf v_0,\dots$ are random variables, so is $\boldsymbol\varphi(k,x,h)$.

We consider in this work the expected infinite-horizon discounted cost over random closed-loop solutions $\boldsymbol\varphi(\cdot,x,h)$ of \eqref{eq:nonlinear_system}.
That is, for any $x\in\mc X$ and $h\in\mc H_a$, let
\begin{align}
    J_{\gamma}(x,h) := \E\left[\sum_{k=0}^\infty \gamma^k\ell(\boldsymbol\varphi(k,x,h), h(\boldsymbol\varphi(k,x,h)))\right]\label{eq:def:costfunction}
\end{align}
be the \textit{cost function}, where $\ell:\mc X\times\mc U\to\R_{\geq0}$ is the \textit{stage cost}, which takes non-negative values, and $\gamma\in(0,1)$ is the \textit{discount factor}.
Note that \eqref{eq:def:costfunction} uses stationary policies, meaning that the same function $h$ is applied at every time step, rather than non-stationary, where the functions may differ between time steps.
Under certain conditions it can be shown that stationary policies achieve the same minimum of the cost function as non-stationary policies \cite[Theorem 4.2.3(d)]{Hern_ndez_Lerma_1996}.

The discount factor $\gamma\in(0,1)$ generally ensures finiteness of $J_\gamma(x,h)$.
In many stochastic systems, in particular with additive noise, no policy leads to convergence to 0 due to persistent random state variation.
Then, the undiscounted ($\gamma=1$) cost function is infinite for any policy, while the discounted ($\gamma\in(0,1)$) cost function is finite for suitable policies keeping expected stage costs below a certain bound.

\subsection{Optimal and near-optimal controls}
For each $\gamma \in (0,1)$ we define the (optimal) value function $V_{\gamma}: \mbox{\rm dom}(V_{\gamma}) \rightarrow \R_{\geq 0}$ with
\begin{align}
    \mbox{\rm dom}(V_{\gamma}):= \left\{ x \in\mc X : \inf_{h \in\mc H_{a}} J_{\gamma}(x,h) < \infty \right\}
\end{align}

and
\begin{align}
    V_{\gamma}(x):=  \inf_{h \in\mc H_{a}} J_{\gamma}(x,h) \label{eq:value_function}
\end{align}
for all $x \in \mbox{\rm dom}(V_{\gamma})$.  
The (near-)optimal feedback policy is based on Bellman equation
\begin{align}
    V_\gamma(x) = \inf_{u\in\mc U}\big(\ell(x,u)+\gamma\E[ V_\gamma(f(x,u,\mathbf v))]\big)~~\forall x\in\mc X,\label{eq:Bellman}
\end{align}
see \cite[Theorem 4.2.3]{Hern_ndez_Lerma_1996} for sufficient conditions for satisfaction of \eqref{eq:Bellman}.
Generally, the infimum may not be attained as a minimum, and in practice only approximate solutions may be found.
To address this phenomenon, our controls can exceed the optimal cost by at most a state-dependent and known error bound $\eta:\mc X\to\R_{\geq0}$ and are hence referred to as \textit{near-optimal}.
With $\ell(x,u)+\gamma\E[V_\gamma(f(x,u,\mathbf v))]$, which is the cost incurred from $u$ followed by optimal control, being desired to be no more than $\eta(x)$ above the minimal cost $V_\gamma(x)$, the near-optimal set-valued feedback policy $\mc U_\gamma^\eta:\mc X\rightrightarrows\mc U$ is defined by
\begin{align}\begin{split}
    \mc U_{\gamma}^\eta(x) := \{u\in\mc U~|~
    &V_{\gamma}(x) + \eta(x)\\
    &\geq \ell(x,u) + \gamma\E[V_{\gamma}(f(x,u,\mathbf v))]\}
\end{split}\label{eq:U_gamma_eta}\end{align}
for any $x\in\mc X$.
For any state $x$, we allow for any control $u\in\mc U_\gamma^\eta(x)$ to be applied to the system, provided a strict causality property is verified, which is introduced next.
Note that \eqref{eq:U_gamma_eta} gives optimal inputs when $\eta(x)=0$ for any $\forall x\in\mc X$, which our results allow.

\subsection{Strictly causal random closed-loop solutions}
We now define solutions to the closed loop of system \eqref{eq:nonlinear_system} and near-optimal feedback \eqref{eq:U_gamma_eta}.
We will study in particular the strictly causal random solutions of such a system, that is, when the control $u\in\mc U_\gamma^\eta(x)$ may depend only on past random inputs, but not on the current and future random inputs.
Not only is this property required to establish our results, but \cite[Example 2]{grammatico_discrete-time_2013} demonstrates that without it robustness to even small perturbations may not be guaranteed.

To formalize strict causality, let the natural filtration $\mc F_{-1}\subseteq \mc F_0\subseteq\dots$ of the random input sequence $(\mathbf v_0, \mathbf v_1,\dots)$ be given by
\begin{align}
    \mc F_k:=\{\{\omega\in\Omega:(\mathbf v_0(\omega),\dots,\mathbf v_k(\omega))\in F\}|F\in\mc B(\R^{q(k+1)})\}
\end{align}
for $k\in\Z_{\geq0}$, where $\mc B(\R^{q(k+1)})$ denotes the Borel $\sigma$-algebra, and $\mc F_{-1}:=\{\varnothing,\Omega\}$.
Every $\mc F_k$ is a sub-$\sigma$-algebra of $\mc F$.
The desired notion of solutions is now defined as follows.
\begin{defn}\label{defn:solution}
    A sequence $\mathbf x = (\mathbf x_0, \mathbf x_1, \dots)$ of random variables $\mathbf x_k:\Omega\to\mc X$ is a \textit{strictly causal random closed-loop solution} of \eqref{eq:nonlinear_system} and \eqref{eq:U_gamma_eta} for initial state $x_0\in\mc X$, denoted by $\mathbf x\in \mc S_\gamma^\eta(x_0)$, if $\mathbf x_0(\omega)=x_0$ for all $\omega\in\Omega$, and for any $k\in\Z_{\geq0}$ there exists an $\mc F_{k-1}$-measurable random variable $\mathbf u_k:\Omega\to\mc U$ such that $\mathbf u_k(\omega)\in\mc U_\gamma^\eta(\mathbf x_k(\omega))$ and $\mathbf x_{k+1}(\omega) = f(\mathbf x_k(\omega),\mathbf u_k(\omega),\mathbf v_k(\omega))$ for all $\omega\in\Omega$.
\ediamond\end{defn}

Note that Definition \ref{defn:solution} implies that $\mathbf x_k$ is $\mc F_{k-1}$-measurable for every $k\in\Z_{\geq0}$.
We use the concept of solutions in Definition \ref{defn:solution} rather than that of Section 3.1 for greater generality, as it allows all control selections in $\mc U_\gamma^\eta$ (in particular, different controls for same states, but at different times), as long as strict causality is maintained.
Nevertheless, any solution $\boldsymbol\varphi(\cdot,x,h)$ for $x\in\mc X$ and where $h\in\mc H_a$ is a selection of $\mc U_\gamma^\eta$ (that is, $h(x)\in\mc U_\gamma^\eta(x)$ for all $x\in\mc X$), qualifies as a solution according to Definition \ref{defn:solution}.
For further reference, see strictly causal generalized random solutions to stochastic difference inclusions in \cite[Section 5]{grammatico_discrete-time_2013}, of which Definition \ref{defn:solution} is a special case.

\section{Semi-global practical recurrence}\label{s:Recurrence}
In this section we present the main result about semi-global practical recurrence, inspired by \cite[Definition 5]{grammatico_discrete-time_2013}, for strictly causal random closed-loop solutions of system \eqref{eq:nonlinear_system} and feedback policy \eqref{eq:U_gamma_eta}, as defined in Definition \ref{defn:solution}, under certain assumptions that are stated in the following.

\subsection{Assumptions}
We start with mild measurability assumptions to ensure that expectations exist.

\begin{assum}\label{assum:measurable}
    The functions $f,\ell$ and $\eta$ are measurable.
    Furthermore, $V_\gamma:\mc X\to\R$ is finite and measurable for any $\gamma\in(0,1)$.
\ediamond\end{assum}

The main assumptions are similar to the ones made in \cite{grimm_model_2005,postoyan_stability_2017} on deterministic systems, but adapted for stochastic systems.
They relate the costs in the optimal control problem to a state space metric $\sigma:\mc X\to\R_{\geq0}$ that recurrence revolves around.
Function $\sigma$ can be defined as $|\cdot|$ or $|\cdot|^2$ when studying recurrence of the origin, or the distance to a subset of $\mc X$ that is desired to be recurrent.
The first condition is cost-controllability, which in \cite{postoyan_stability_2017} guarantees existence of a $\mc K_\infty$-function of $\sigma(x)$ that upper bounds the value function $V_\gamma(x)$ for all discount factors $\gamma\in(0,1)$.
The second condition is detectability of the stage cost, and, roughly speaking, means, that small stage costs imply that $\sigma(x)$ becomes small eventually over time along solutions to \eqref{eq:nonlinear_system}.
After stating the assumptions now, we explain them in more detail and illustrate them on an example below.

\begin{assum}\label{assm:cost_controllability_detectability}
    Let $\sigma:\mc X\to\R_{\geq0}$ be measurable, locally bounded and radially unbounded.
    \begin{enumerate}[label=(\roman*)]
        \item \textit{Cost-controllability with lower bound:} There exist $\overline\alpha_V\in\mc K_\infty$ and $e_1,e_2\geq0$ such that for any $\gamma\in(0,1)$ there exists $c_\gamma\geq0$ such that for any $x\in\mc X$,
        \begin{align}-e_1\leq V_\gamma(x) - \frac{c_\gamma}{1-\gamma}\leq e_2 + \overline\alpha_V(\sigma(x)).\end{align}
        Furthermore, there exists $c\geq0$ such that $c_\gamma\leq c$ for any $\gamma\in(0,1)$.\label{item:assm:cost_controllability}
        \item \textit{Detectability:} There exist a measurable function $W:\mc X\to\R_{\geq0}$, functions $\overline\alpha_W,\alpha_W\in\mc K_\infty$ and $d>0$, such that for any $x\in\mc X$ and $u\in\mc U$,
        \begin{align}\begin{split}
        W(x) &\leq \overline\alpha_W(\sigma(x)),\\
        \E[W(f(x,u,\mathbf v))] - W(x) &\leq -\alpha_W(\sigma(x)) + \ell(x,u) + d.\label{eq:dissipation}
        \end{split}\end{align}\label{item:assm:detectability}
        \ediamond
    \end{enumerate}
\end{assum}

Assumptions \ref{assm:cost_controllability_detectability} generalizes \cite[Assumption 1]{postoyan_stability_2017} from deterministic to stochastic systems through the parameters $c,d,e_1,e_2\geq0$.
These parameters are designed to capture types of behaviour unique to stochastic systems, and cover a broad range of ways the stochastic input impacts the system, like additive and multiplicative noise.
The parameters $c$ and $d$, as well as the near-optimality bound $\eta$, will appear in the statements of our recurrence results, in the way that lower values yield a smaller recurrent set.
When $c, d, e_1, e_2$ and $\eta$ are all zero, \cite[Assumption 1]{postoyan_stability_2017} is recovered.

Now let us explain Assumption \ref{assm:cost_controllability_detectability} in detail.
Regarding item \ref{item:assm:cost_controllability}, for $x\in\mc X$ with $\sigma(x)=0$ the value function $V_\gamma(x)$ differs from $c_\gamma/(1-\gamma)$ only by the constants $e_1,e_2$, and for general $x\in\mc X$ the upper bound may increase through $\overline\alpha_V$ as $\sigma(x)$ increases.
An equivalent formulation of part \eqref{assm:cost_controllability_detectability} is obtained if $c_\gamma/(1-\gamma)$ and $e_1$ are replaced by $\inf_{x\in\mc X}V_\gamma(x)$ and $0$, respectively.
However, $c_\gamma$ will appear in the main theorem and the presented version with $e_1$ allows for more choices for $c_\gamma$ and easier verification.
The reason to include $c_\gamma$ is the same given in Section \ref{s:3.1:Classofsystemsandcosts} for discounting.
Random inputs, particularly additive noise, can lead to a certain minimum expected stage cost at every time instant, causing the value function $V_\gamma$ to blow up as $\gamma$ approaches 1.
This is captured in the $c_\gamma/(1-\gamma)$ term, where $c_\gamma$ can be interpreted as a minimum asymptotic average stage cost and is required to be bounded for $\gamma\in(0,1)$.
The division by $1-\gamma$ comes from the geometric series, when these costs $c_\gamma$ are added up with discount factor $\gamma$.
Note that $e_1,e_2$ and $\overline\alpha_V$ need to be independent of $\gamma$.
This means that, although $V_\gamma(x)$ may blow up as $\gamma\to1$, the differences $V_\gamma(x_1)-V_\gamma(x_2)$ for same $\gamma$ and different $x_1,x_2\in\mc X$ cannot.
The example below illustrates that this is a reasonable assumption.





Moving on to item \ref{item:assm:detectability} of Assumption \ref{assm:cost_controllability_detectability}, $W$ can be chosen constantly zero if there exists $\alpha_W\in\mc K_\infty$ such that $\ell(x,u)\geq\alpha_W(\sigma(x))$ for all $(x,u)\in\mc X\times\mc U$.
The formulation with $W$ is more general as it allows positive semi-definite $\ell$ that may depend on only some of the state components, for example.
Inequality \eqref{eq:dissipation} resembles a dissipation-type inequality with storage function $W$ and supply rate $\ell$.
While the special case $W=0$ provides that $\sigma(x)$ is small right away when $\ell(x,u)$ is small, in the general case small $\ell(x,u)$ and large $\sigma(x)$ mean that $W$ decreases in expectation over time due to \eqref{eq:dissipation}, but since it cannot become negative, $\sigma(x)$ must become small eventually.
Keeping to the analogy of dissipativity, one may expect the random input to supply some average amount of energy at every time instant, which is captured in $d$.

To see how Assumption \ref{assm:cost_controllability_detectability} is consistent with and generalizes the linear stochastic case studied in \cite{granzotto_stability_2024}, consider the next example.

\begin{exmp}\label{exmp:LQR}
    Consider the linear system
    \begin{align}
        x^+ &= f(x,u,v) := Ax + Bu + Lv
    \end{align}
    with matrices $A\in\R^{n\times n}, B\in\R^{n\times m}, L\in\R^{n\times q}$, and the stage cost
    \begin{align}
        \ell(x,u) = x^\top Q x + u^\top R u,\label{eq:exmp:stagecost}
    \end{align}
    where $Q\in\R^{n\times n}$ is positive semi-definite and $R\in\R^{m\times m}$ is positive definite.
    We are interested in recurrence of (open neighbourhoods of) the origin, and therefore choose $\sigma(x):=|x|$.
    If $(A,B)$ is stabilizable and $(Q^{1/2},A)$ is detectable then, for any $\gamma\in(0,1)$, the algebraic Riccati equation
    \begin{align*}
        P_\gamma = Q + \gamma A^\top P_\gamma A - \gamma^2A^\top P_\gamma B(R+\gamma B^\top P_\gamma B)^{-1}B^\top P_\gamma A
    \end{align*}
    has a positive definite solution $P_\gamma\in\R^{n\times n}$, and the value function takes the form
    \begin{align}
        V_\gamma(x) = x^\top P_\gamma x + \frac\gamma{1-\gamma}\E[\mathbf v^\top L^\top P_\gamma L\mathbf v]
    \end{align}
    for any $x\in\R^n$ \cite[Chapter 6.2]{davis_stochastic_1985}. If $P_\gamma$ is bounded over $\gamma\in(0,1)$, then item \ref{item:assm:cost_controllability} of Assumption \ref{assm:cost_controllability_detectability} is verified with $e_1:=e_2:=0$, 
    \begin{align}
        c_\gamma:=\gamma\E[\mathbf v^\top L^\top P_\gamma L\mathbf v], \quad\overline\alpha_V(s):=\sup_{\gamma\in(0,1)}||P_\gamma||_2s^2
    \end{align}
    (recall that $||\cdot||_2$ denotes the spectral norm, i.e., the largest singular value, of a matrix).
    Note that, despite $V_\gamma(x)\to\infty$ as $\gamma\to1$, the difference $V_\gamma(x_1)-V_\gamma(x_2)=x_1^\top P_\gamma x_1-x_2^\top P_\gamma x_2$ depends only on the quadratic term and therefore remain bounded as $\gamma\to1$.
    Regarding detectability, if $(Q^{1/2},A)$ is detectable, by 
    \cite[Appendix B]{postoyan_stability_2017}
    there exist a matrix $M\in\R^{n\times n}$ and a constant $a_W\in\R_{>0}$ such that item \ref{item:assm:detectability} of Assumption \ref{assm:cost_controllability_detectability} holds with $W(x) := x^\top M x, \alpha_W(s) := a_Ws^2$, and $d:=\E[\mathbf v^\top L^\top ML\mathbf v]$.
\ediamond\end{exmp}

Finally, we make the next mild assumption.
\begin{assum}\label{assm:locally_bounded} The function $f$ and the set-valued mapping
        \begin{align}
            \mc U_{(0,1)}^\eta:\mc X\rightrightarrows\mc U,\quad x\mapsto \bigcup_{\gamma\in(0,1)}\mc U_\gamma^\eta(x) =: \mc U_{(0,1)}^\eta(x) \label{eq:assm:U(0,1)definition}
        \end{align}
        are locally bounded, 
        which for the set-valued mapping $\mc U_{(0,1)}^\eta$ means that,
        for any compact set $D\subseteq\mc X$, the image $\mc U_{(0,1)}^\eta(D):=\bigcup_{x\in D}\mc U_{(0,1)}^\eta(x)$ is bounded.
\ediamond\end{assum}

To ease the test of Assumption \ref{assm:locally_bounded}, the next proposition gives sufficient conditions for local boundedness of $\mc U_{(0,1)}^\eta$.
\begin{prop}\label{prop:U01_locally_bounded}
    Let Assumption \ref{assum:measurable} and \ref{assm:cost_controllability_detectability} hold. Further assume that $\eta$ is locally bounded and that $\ell$ is level-bounded in $u$ locally uniformly in $x$, i.e., for each compact set $D\subseteq\mc X$ the function $\ell_D:\mc U\to\R_{\geq0}, u\mapsto\inf_{x\in D}\ell(x,u)=:\ell_D(u)$ is radially unbounded. 
    Then, the set-valued mapping $\mc U_{(0,1)}^\eta$ defined in \eqref{eq:assm:U(0,1)definition} is locally bounded.
\ediamond\end{prop}
The level-boundedness condition on $\ell$ (also see \cite[Definition\,1.16]{RockWets98}) penalizes large control effort and is fulfilled for the quadratic costs \eqref{eq:exmp:stagecost} due to $R$ being positive definite.
The proof of Proposition \ref{prop:U01_locally_bounded} can be found in the appendix.

\subsection{Main result}
The main result of this section is presented in two parts: a lemma and the main recurrence theorem. These two parts are adaptions of items (i) and (ii) of the definition of semi-global practical recurrence in \cite[Definition 5]{grammatico_discrete-time_2013}, respectively.
The lemma states that for any $\Delta_0>0$, there exists a bounded set of the form $\sigma^{-1}(\Delta\B)$ for some $\Delta>0$, in which closed-loop solutions starting in $\sigma^{-1}(\Delta_0\B)$ stay for a desired time $\TIMETEMP$ with at least a desired probability $1-p$.
Hereby, $\Delta$ is independent of $\gamma$.
This is formalized below.

\begin{lem}\label{lem:Boundedness_with_high_prob}
    Let Assumptions \ref{assum:measurable}, \ref{assm:cost_controllability_detectability} and \ref{assm:locally_bounded} hold. Then,
    \begin{equation}
    \addtolength{\fboxsep}{5pt}
    \begin{gathered}\forall\Delta_0>0~\forall p>0~\forall \TIMETEMP\in\Z_{\geq0}~\exists\Delta>0\\
    \forall\gamma\in(0,1)~\forall x_0\in\sigma^{-1}(\Delta_0\B)~\forall \mathbf x\in\mc S_\gamma^\eta(x_0):\\
    \P[\sigma(\mathbf x_k)<\Delta~\forall k\in\{0,\dots,\TIMETEMP\}]>1-p.
    \end{gathered}
    \end{equation}
    
\vspace{-.62cm}\ediamond\end{lem}

The proof can be found in the appendix.

Leading up to the main theorem, for $\gamma\in(0,1)$ and $\delta>0$ we define the set
\begin{align}\mc A_{\gamma,\delta}:= \big\{x\in\mc X~|~\alpha_W(\sigma(x))\leq c_\gamma+d+\eta(x)+\delta\big\}.\end{align}
We will show semi-global practical recurrence of these sets, in the sense that for any $p>0,\delta>0$ and any maximum value $\Delta_0$ of $\sigma(x_0)$ for potential initial states $x_0$, it can be ensured with $\gamma$ close enough to 1 that the closed-loop solutions $\mathbf x\in\mc S_\gamma^\eta(x_0)$ enter $\mc A_{\gamma,\delta}$ with probability at least $1- p$ within some time $\TIMETEMP$ independent of $\gamma,x_0$ and $\mathbf x$.
``Semi-global" refers to the fact that, while for fixed $\gamma$ there can be $x_0$ (that have high $\sigma(x_0)$) for which recurrence does not hold, for any fixed $x_0$ recurrence does hold if $\gamma$ is close enough to 1.
``Practical" is the analogous concept for approaching $\mc A_{\gamma,0}$.
While for fixed $\gamma$ the recurrent set generally needs to be larger than $\mc A_{\gamma,0}$, namely $\mc A_{\gamma,\delta}$, $\delta$ can be arbitrarily small if $\gamma$ is close enough to 1.

\begin{thm}\label{thm:Recurrence}
    Let Assumptions \ref{assum:measurable}, \ref{assm:cost_controllability_detectability} and \ref{assm:locally_bounded} hold. Then,
   \begin{equation}
    \begin{gathered}\forall\Delta_0,\delta>0~\forall p>0~\exists \TIMETEMP\in\Z_{\geq0}~\exists\gamma^\star\in(0,1)\\ \forall \gamma\in(\gamma^\star,1)~\forall x_0\in\sigma^{-1}(\Delta_0\B)~\forall \mathbf x\in\mc S_\gamma^\eta(x_0):\\
        \P[\exists k\in\{0,\dots,\TIMETEMP\}:\mathbf x_k\in\mc A_{\gamma,\delta}]>1- p.
    \end{gathered}
    \end{equation}
    
\vspace{-.62cm}\ediamond\end{thm}

The proof of Theorem \ref{thm:Recurrence} can be found in the appendix.
Note that, since the system \eqref{eq:nonlinear_system} is time-invariant, Lemma \ref{lem:Boundedness_with_high_prob} and Theorem \ref{thm:Recurrence} can be combined to guarantee recurrence farther into the future as well, as for any time $\TIMETEMP_0\in\Z_{\geq0}$, Lemma \ref{lem:Boundedness_with_high_prob} can be used to construct a set in which solutions stay until $\TIMETEMP_0$ with high probability, starting from which Theorem \ref{thm:Recurrence} can guarantee recurrence in overall time within $\{\TIMETEMP_0,\dots,\TIMETEMP_0+\TIMETEMP\}$ with high probability.

Next, Theorem \ref{thm:Recurrence} is generalized to the case with small strictly causal perturbations.

\section{Robustness}\label{s:Robust_Recurrence}
In this section, we state robust semi-global practical recurrence under additional continuity assumptions.
Solutions are defined similarly to Definition \ref{defn:solution}, except that in addition to $\mathbf u_k$, perturbation random variables $\mathbf p_k^i, i\in\{1,2,3,4\}$ are introduced.
\begin{defn}\label{defn:epsilon_perturbed_solution}

    For $\varepsilon>0$, a sequence $\mathbf x = (\mathbf x_0, \mathbf x_1, \dots)$ of random variables $\mathbf x_k:\Omega\to\mc X$ is an  \textit{$\varepsilon$-perturbed strictly causal random closed-loop solution} of \eqref{eq:nonlinear_system} and \eqref{eq:U_gamma_eta} for initial state $x_0\in\mc X$, denoted by $\mathbf x\in \mc S_{\gamma,\varepsilon}^\eta(x_0)$, if $\mathbf x_0=x_0$, and for any $k\in\Z_{\geq0}$ there exist $\mc F_{k-1}$-measurable $\mathbf u_k:\Omega\to\mc U, \mathbf p_k^2:\Omega\to\varepsilon\B^n, \mathbf p_k^3:\Omega\to\varepsilon\B^m$ 
    such that
    \begin{align}
        \mathbf x_k+\mathbf p_k^2\in\mc X,\quad\mathbf u_k\in\mc U_\gamma^\eta(\mathbf x_k+\mathbf p_k^2)+\mathbf p_k^3,~~~~~~~~~~~
    \end{align}
    and $\mc F_k$-measurable $\mathbf p_k^1,\mathbf p_k^4:\Omega\to\varepsilon\B^n$ such that
    \begin{align}
        \mathbf x_k+\mathbf p_k^1\in\mc X,\quad\mathbf x_{k+1} = f(\mathbf x_k+\mathbf p_k^1,\mathbf u_k,\mathbf v_k)+\mathbf p_k^4.
    \end{align}
\ediamond\end{defn}
Perturbations $\mathbf p_k^2$ and $\mathbf p_k^3$ apply to the input and output of the controller, and have to be $\mc F_{k-1}$-measurable like $\mathbf u_k$, meaning that they can depend only on $\mathbf v_0,\dots,\mathbf v_{k-1}$.
On the other hand, $\mathbf p_k^1$ and $\mathbf p_k^4$ apply on the state before and after $f$ is applied, and need only to be $\mc F_k$-measurable, meaning they can depend on $\mathbf v_k$ as well.

We make the following assumption strengthening the measurability of Assumption \ref{assum:measurable} to continuity.
\begin{assum}\label{assm:continuity}
    The following holds.
    \begin{enumerate}[label=(\roman*)]
        \item $\ell,\sigma$ and $\eta$ are continuous.\label{item:assm:robust1}\label{item:assm:robust2}\label{item:assm:robust2p5}
        \item $(x,u)\mapsto f(x,u,v)$ is continuous for any $v\in\mc V$.\label{item:assm:robust3}
        \item $\{V_\gamma~|~\gamma\in(0,1)\}$ and $\{\E[V_\gamma(f(\cdot,\cdot,\mathbf v))]~|~\gamma\in(0,1)\}$ are pointwise equi-continuous, i.e., $\forall\varepsilon>0~\forall x\in\mc X~\exists\delta>0~\forall y\in\mc X$,
        \begin{align}
        \begin{gathered}
            |x-y|<\delta\Ra \forall\gamma\in(0,1):|V_\gamma(x)-V_\gamma(y)|<\varepsilon,
        \end{gathered}
        \end{align}
        similar for $\{\E[V_\gamma(f(\cdot,\cdot,\mathbf v))]~|~\gamma\in(0,1)\}$. Furthermore, $W$ is continuous.\label{item:assm:robust4}
        \ediamond
    \end{enumerate}
\end{assum}
All conditions of Assumption \ref{assm:continuity} are satisfied in the LQR problem of Example \ref{exmp:LQR}.


We can now state the robust versions of Lemma \ref{lem:Boundedness_with_high_prob} and Theorem \ref{thm:Recurrence}.
The difference in the next results is in the existence of a perturbation bound $\varepsilon$ (that generally depends on $\Delta_0,\delta, p$ and $\TIMETEMP$) such that desired properties are fulfilled even for $\varepsilon$-perturbed, not just exact, strictly causal random closed-loop solutions.
Parameter $\varepsilon$ is independent of $\gamma$ due to the pointwise equi-continuity in Assumption \ref{assm:continuity}.

\begin{lem}\label{lem:r:Boundedness_with_high_prob}
    Let Assumptions \ref{assum:measurable}, \ref{assm:cost_controllability_detectability}, \ref{assm:locally_bounded} and \ref{assm:continuity} hold. Then,
    \begin{equation}
    \addtolength{\fboxsep}{5pt}
    \begin{gathered}\forall\Delta_0>0~\forall p>0~\forall \TIMETEMP\in\Z_{\geq0}~\exists\Delta>0~\exists\varepsilon>0\\
    \forall\gamma\in(0,1)~\forall x_0\in\sigma^{-1}(\Delta_0\B)~\forall \mathbf x\in\mc S_{\gamma,\varepsilon}^\eta(x_0):\\
    \P[\sigma(\mathbf x_k)<\Delta~\forall k\in\{0,\dots,\TIMETEMP\}]>1-p.
    \end{gathered}
    \end{equation}
    
\vspace{-.62cm}\ediamond\end{lem}


\begin{thm}\label{thm:r:Recurrence}
    Let Assumptions \ref{assum:measurable}, \ref{assm:cost_controllability_detectability}, \ref{assm:locally_bounded} and \ref{assm:continuity} hold. Then,
    \begin{equation}
    \begin{gathered}\forall\Delta_0,\delta>0~\forall p>0~\exists \TIMETEMP\in\Z_{\geq0}~\exists\varepsilon>0~\exists\gamma^\star\in(0,1)\\
    \forall \gamma\in(\gamma^\star,1)~\forall x_0\in\sigma^{-1}(\Delta_0\B)~\forall \mathbf x\in\mc S_{\gamma,\varepsilon}^\eta(x_0):\\
        \P[\exists k\in\{0,\dots,\TIMETEMP\}:\mathbf x_k\in\mc A_{\gamma,\delta}]>1-p.
    \end{gathered}
    \end{equation}

\vspace{-.62cm}\ediamond\end{thm}
The proofs of Lemma \ref{lem:r:Boundedness_with_high_prob} and Theorem \ref{thm:r:Recurrence} follow a similar structure to those of Lemma \ref{lem:Boundedness_with_high_prob} and Theorem \ref{thm:Recurrence}.

\section{Conclusion}\label{s:Conclusion}
For discrete-time stochastic systems controlled by infinite-horizon discounted optimal control inputs, appropriate cost-controllability and detectability assumptions imply semi-global practical recurrence with respect to the discount factor.
This property is uniform in the sense that the time to enter the recurrent set with high probability is independent of discount factor and initial state.
Under additional continuity assumptions, robustness is guaranteed.
Ongoing work is being done to obtain recurrence properties under relaxations of Assumption \ref{assm:cost_controllability_detectability}, as well as for the undiscounted case and stochastic model predictive control.

\bibliography{references}

\appendix
\section{Proof of Proposition \ref{prop:U01_locally_bounded}}
Let $D\subseteq\mc X$ be compact.
We show that $\mc U_{(0,1)}^\eta(D)$ as defined in Assumption \ref{assm:locally_bounded} is bounded.
Because $D$ is compact and $\sigma$, $\overline\alpha_V$ and $\eta$ are locally bounded, there exists $R\in\R$ such that $c+e_1+e_2+\overline\alpha_V(\sigma(x'))+\eta(x')< R$ for any $x'\in D$.

Suppose $u\in\mc U_{(0,1)}^\eta(D)$, then by \eqref{eq:assm:U(0,1)definition} there exist $x\in D$ and $\gamma\in(0,1)$ such that $u\in\mc U_\gamma^\eta(x)$.
By applying \eqref{eq:U_gamma_eta} as well as part \ref{item:assm:cost_controllability} of Assumption \ref{assm:cost_controllability_detectability},
\begin{align*}
    \ell_D(u)&\leq\ell(x,u) \leq V_\gamma(x)+\eta(x) - \gamma\E V_\gamma(f(x,u,\mathbf v))\\
    &\leq \frac{c_\gamma}{1-\gamma}+e_2+\overline\alpha_V(\sigma(x))+\eta(x)-\frac{\gamma c_\gamma}{1-\gamma}+\gamma e_1\\
    &\leq c+e_1+e_2 + \overline\alpha_V(\sigma(x))+\eta(x)<R.
\end{align*}

We have proved $\ell_D(u)<R$ for any $u\in\mc U_{(0,1)}^\eta(D)$, which implies $\mc U_{(0,1)}^\eta(D)\subseteq\ell_D^{-1}([0,R))$.
Because $\ell_D$ is radially unbounded, $\ell_D^{-1}([0,R))$ and hence $\mc U_{(0,1)}^\eta(D)$ are bounded.
This completes the proof.


\section{Proof of Lemma \ref{lem:Boundedness_with_high_prob}}\label{s:proof:lem:Boundedness}
    The proof is similar to that of \cite[Lemma 4]{teel_converse_2014}.
    Given $\Delta_0>0, \TIMETEMP\in\Z_{\geq0}$ and $p>0$, let $\tmc V\subseteq\R^q$ be compact such that $\P[\mathbf v\notin\tmc V]< p/\TIMETEMP$.
    We define sets $\mc X_0,\dots,\mc X_\TIMETEMP$ as
    \begin{align*}
        \mc X_0 &:= \cl\left(\sigma^{-1}(\Delta_0\B)\right),\\
        \mc X_{i+1} &:= \cl\left(f\left(\mc X_i, \cl(\mc U_{(0,1)}^\eta(\mc X_i)),\tmc V\right)\right)~\forall i\in\{0,\dots,\TIMETEMP-1\}.
    \end{align*}
    Note that $\mc X_i\subseteq\mc X$ for any $i\in\{0,\dots,T\}$ because $f$ maps to $\mc X$, and $\mc X$ and $\mc U$ are closed.
    We assert that each of these sets is compact. In fact, $\mc X_0$ is compact because $\sigma$ is radially unbounded.
    Now suppose that $\mc X_i$ is compact, then, because $\mc U_{(0,1)}^\eta$ is locally bounded and $\mc U$ is closed, $\cl(\mc U_{(0,1)}^\eta(\mc X_i))\subseteq\mc U$ is compact.
    Using this as well as compactness of $\mc X_i$ and $\tmc V$, it follows from local boundedness of $f$ that $f\left(\mc X_i, \cl(\mc U_{(0,1)}^\eta(\mc X_i)),\tmc V\right)$ is bounded, and $\mc X_{i+1}$ is compact.

    Choose $\Delta>\sup_{x\in\mc X_0\cup\dots\cup\mc X_\TIMETEMP}\sigma(x)$ (this is possible because $\mc X_0,\dots,\mc X_\TIMETEMP$ are compact and $\sigma$ is locally bounded).
    Let $x_0\in\Delta_0\B$ and $\mathbf x\in\mc S_\gamma^\eta(x_0)$.
    Then, for any $\omega\in\Omega$, $\mathbf v_k(\omega)\in\tmc V$ for any $k\in\{0,\dots,\TIMETEMP-1\}$ implies $\mathbf x_k(\omega)\in\mc X_k$ for any $k\in\{0,\dots,\TIMETEMP\}$, and consequently $\mathbf \sigma(\mathbf x_k(\omega))<\Delta$ for any $k\in\{0,\dots,\TIMETEMP\}$.
    Hence,
    \begin{align*}
        &\P[\sigma(\mathbf x_k)<\Delta~\forall k\in\{0,\dots,\TIMETEMP\}] \\&\geq \P[\mathbf v_k\in\tmc V~\forall k\in\{0,\dots,\TIMETEMP-1\}]\\
        &\geq 1-\P[\mathbf v_0\notin\tmc V] - \dots - \P[\mathbf v_{\TIMETEMP-1}\notin\tmc V]\\
        &>1-p,
    \end{align*}
    where we used the law of total probability.

\omitcustom{
    The reader may skip ahead to Parts 2 and 4, which contain the main ideas of the proof, and follow the references therein to Parts 1 and 3 for derivations of the bounds used in Part 4.\\
    \underline{Part 1: Setup.}
    Define $ Y_\gamma := V_\gamma+\frac1\gamma W-\frac{c_\gamma}{1-\gamma}+c_\gamma+d+e_1$, $e:=e_1+e_2$ and $\overline\alpha_Y:=\overline\alpha_V+2\overline\alpha_W$ as in Lemma \ref{lem:Sandwich_bounds}.
    Given $\Delta_0>0,  p>0$, $\TIMETEMP\in\Z_{\geq0}$ and $\varepsilon>0$, choose 
    \begin{align}
    \Delta := \max\left\{\frac{\frac1{1- p}(\overline\alpha_Y(\Delta_0)+c+d+e+\TIMETEMP(2(\bar\eta+c+d)+1))}{1-(1- p)^{1/\TIMETEMP}}, \Delta_0\right\} \label{Lemmadag1}\tag{$\dagger_1$}
    \end{align}
    and $\gamma^\star\in(1/2,1)$ such that
    \begin{align}
        \frac1{\gamma^\star}(1-\gamma^\star)(\overline\alpha_V(\Delta)+e_2)<1. \label{Lemmadag2}\tag{$\dagger_2$}
    \end{align}
    Let $\gamma\in(\gamma^\star,1)$, $x_0\in\sigma^{-1}(\Delta_0\B)$ and $(\mathbf x_k,\mathbf u_k)$ be a solution of \eqref{eq:closed_loop_difference_inclusion} with discount factor $\gamma$ and initial state $x_0$.
    For any $k\in\Z_{\geq0}$ define
    \begin{align}\mc R_k := \{\omega\in\Omega~|~\forall j=0,\dots,k:\sigma(\mathbf x_j(\omega))<\Delta\}.\end{align}
    Because $\sigma$ and $\mathbf x_j$ are measurable, these sets are measurable.
    Note that $\P(\mc R_0)=1$ since $\Delta\geq\Delta_0$.
        
    \noindent\underline{Part 2: Assertion for proof by induction.} 
    We will show that for any $k\in\{0,\dots,\TIMETEMP-1\}:$
    \begin{align}\frac{\P(\mc R_{k+1})}{\P(\mc R_k)} > (1- p)^{1/\TIMETEMP}.\end{align}
    Note that this induction proof will also show that $\P(\mc R_k)>0$ for any $k\in\{0,\dots,\TIMETEMP\}$, which means that conditional expectations on these events exist.

    \noindent\underline{Part 3: Derivation of bounds along the closed loop.}
    Let $k\in\Z_{\geq0}$ such that $\P(\mc R_k)>0$.
    \underline{3.1: Upper bound for $\E[ Y_\gamma(\mathbf x_{k+1})|\mc R_k]$.} In this part, we derive a bound for the growth of the expected value function, conditioned on $\mc R_k$. 
    Due to strict causality, $\mathbf x_k$ and $\mathbf u_k$ are $\mc F_{k-1}$-measurable (where $\mc F_{-1}:=\{\emptyset,\Omega\}$). Since $\mathbf v_j, j\in\Z_{\geq0}$ are i.i.d., $\mathbf v_k$ is independent of $\mc F_{k-1}$, and thus independent of $\mathbf x_k$ and $\mathbf u_k$.
    As $\mc R_k$ is $\mc F_{k-1}$-measurable as well, $\mathbf v_k$ is also independent of $\mc R_k$.
    Thus
    \begin{align}
        \P^{\mathbf x_k,\mathbf u_k,\mathbf v_k|\mc R_k} = \P^{\mathbf v_k}\otimes\P^{\mathbf x_k,\mathbf u_k|\mc R_k},
    \end{align}
    where $\otimes$ denotes the product measure.
    Using this along with Fubini's theorem yields
    \begin{align}
        \E[ Y_\gamma(\mathbf x_{k+1})|\mc R_k] - \E[ Y_\gamma(\mathbf x_k)|\mc R_k] &= \int \E Y_\gamma(f(x,u,\mathbf v_k)) -  Y_\gamma(x)\td\P^{\mathbf x_k,\mathbf u_k|\mc R_k}(x,u)\\
        &=\int \E V_\gamma(f(x,u,\mathbf v_k))-V_\gamma(x) + \frac1\gamma(\E W(f(x,u,\mathbf v_k))-W(x))\td\P^{\mathbf x_k,\mathbf u_k|\mc R_k}(x,u)
    \end{align}
    By \eqref{eq:U_gamma_epsilon} and \eqref{eq:closed_loop_difference_inclusion}, $V_\gamma(x) + \eta(x)\geq \ell(x,u) + \gamma\E V_\gamma(f(x,u,\mathbf v_k))$ for $\P^{\mathbf x_k,\mathbf u_k}$-a.e. $(x,u)$.
    Using this as well as part \ref{item:assm:detectability} of Assumption \ref{assm:cost_controllability_detectability},
    \begin{align}
        \E[ Y_\gamma(\mathbf x_{k+1})|\mc R_k] - \E[ Y_\gamma(\mathbf x_k)|\mc R_k] &\leq \int \frac1\gamma((1-\gamma)V_\gamma(x) +\underbrace{\eta(x) - \alpha_W(\sigma(x))}_{\leq\bar\eta} + d)\td\P^{\mathbf x_k,\mathbf u_k|\mc R_k}(x,u)\\
        &= \int\frac1\gamma((1-\gamma)(\overline\alpha_V(\Delta)+e_2)+\bar\eta+c_\gamma+d)\td\P^{\mathbf x_k,\mathbf u_k|\mc R_k}(x,u)\\
        &= \frac1\gamma((1-\gamma)(\overline\alpha_V(\Delta)+e_2)+\bar\eta+c_\gamma+d),
    \end{align}
    where we used \ref{item:assm:cost_controllability} and \ref{item:assm:eta_growth} of Assumption \ref{assm:cost_controllability_detectability} and $\sigma(x)<\Delta$ for $\P^{\mathbf x_k,\mathbf u_k|\mc R_k}$-a.e. $(x,u)$.
    With \eqref{Lemmadag2}, $\gamma>\gamma^\star$, $\gamma\geq1/2$ and $c_\gamma\leq c$, it follows that
    \begin{align}\E[ Y_\gamma(\mathbf x_{k+1})|\mc R_k] - \E[ Y_\gamma(\mathbf x_k)|\mc R_k]<2(\bar\eta+c+d)+1.\label{Lemmastar1}\tag{$\star_1$}
    \end{align}
    \underline{3.2: Lower bound for $\P(\mc R_{k+1})$.}
    Using Markov's inequality, as well as the left inequality in Lemma \ref{lem:Sandwich_bounds}, we can derive a bound for the probability that $\sigma(\mathbf x_{k+1})>\Delta$ in terms of the expression we just bounded in \eqref{Lemmastar1}:
    \begin{align}
        1-\frac{\P(\mc R_{k+1})}{\P(\mc R_k)} &= \P[\Omega\setminus\mc R_{k+1}|\mc R_k]\\
        &= \P[\sigma(\mathbf x_{k+1})\geq\Delta|\mc R_k]\\
        &\leq \P[ Y_\gamma(\mathbf x_{k+1})\geq\alpha_W(\Delta)|\mc R_k]\\&\leq\frac{\E[ Y_\gamma(\mathbf x_{k+1})|\mc R_k]}{\alpha_W(\Delta)}. \label{Lemmastar2}\tag{$\star_2$}
    \end{align}
    \underline{3.3: Upper bound for $\E[ Y_\gamma(\mathbf x_{k+1})|\mc R_{k+1}]$.}
    Furthermore, since $\mc R_{k+1}\subseteq\mc R_k$ and $ Y_\gamma\geq0$ by Lemma \ref{lem:Sandwich_bounds}, if $\P(\mc R_{k+1})>0$,
    \begin{align}
        \E[ Y_\gamma(\mathbf x_{k+1})|\mc R_{k+1}] &= \P(\mc R_{k+1})^{-1}\int_{\mc R_{k+1}} Y_\gamma(\mathbf x_{k+1})\td\P\\
        &\leq \P(\mc R_{k+1})^{-1}\int_{\mc R_{k}} Y_\gamma(\mathbf x_{k+1})\td\P\\
        &= \frac{\P(\mc R_{k})}{\P(\mc R_{k+1})}\E[ Y_\gamma(\mathbf x_{k+1})|\mc R_k]. \label{Lemmastar3}\tag{$\star_3$}
    \end{align}
    This means that, if $\P(\mc R_{k+1})$ is high (which we can guarantee with \eqref{Lemmastar1} and \eqref{Lemmastar2}), $\E[ Y_\gamma(\mathbf x_{k+1})|\mc R_{k+1}]$ is not much larger than $\E[ Y_\gamma(\mathbf x_{k+1})|\mc R_{k}]$ (which goes into \eqref{Lemmastar1} for the next timestep, and so on).
    
    \noindent\underline{3.4: Upper bound for $\E[ Y_\gamma(\mathbf x_0)|\mc R_0]$.}
    Finally, we need an initial bound.
    Since $\sigma(x_0)<\Delta_0$, it follows from the upper bound in Lemma \ref{lem:Sandwich_bounds} that \begin{align}\E[ Y_\gamma(\mathbf x_0)|\mc R_0] =  Y_\gamma(x_0) < \overline\alpha_Y(\Delta_0)+c+d+e.\end{align}
    \noindent\underline{3.7: Multiple time steps.} With this and applying (\ref{Lemmastar1}) and (\ref{Lemmastar3}) repeatedly, we obtain for $k<\TIMETEMP$,
    \begin{align}
        \E[ Y_\gamma(\mathbf x_{k+1})|\mc R_k] &\leq \frac{\P(\mc R_0)}{\P(\mc R_k)}(\overline\alpha_Y(\Delta_0)+c+d+e) + \underbrace{\left(\frac{\P(\mc R_0)}{\P(\mc R_k)} + \frac{\P(\mc R_1)}{\P(\mc R_k)} + \dots + \frac{\P(\mc R_{k})}{\P(\mc R_k)}\right)}_{\leq (k+1)\frac{\P(\mc R_0)}{\P(\mc R_k)}\leq \TIMETEMP\frac{\P(\mc R_0)}{\P(\mc R_k)}}(2(\bar\eta+c+d)+1)\\
        &\leq \frac{\P(\mc R_0)}{\P(\mc R_k)}(\overline\alpha_Y(\Delta_0)+c+d+e+\TIMETEMP(2(\bar\eta+c+d)+1)).\label{Lemmastar4}\tag{$\star_4$}
    \end{align}
\mbox{}\hfill$\Box$    \underline{Part 4: Induction proof.} Let $k\in\{0,\dots,\TIMETEMP-1\}$ and assume
    \begin{align}\frac{\P(\mc R_{j+1})}{\P(\mc R_j)}>(1- p)^{1/\TIMETEMP} \text{ for every }j\in\{0,\dots,k-1\}.\end{align}
    Then $\P(\mc R_k)>0$, and thus the inequalities from Part 3 can be applied. Combining the induction hypothesis with (\ref{Lemmastar4}) yields
    \begin{align}\E[ Y_\gamma(\mathbf x_{k+1})|\mc R_k] &\leq (1- p)^{-k/\TIMETEMP}(\overline\alpha_Y(\Delta_0)+c+d+e+\TIMETEMP(2(\bar\eta+c+d)+1))\\
    &< \frac1{1- p}(\overline\alpha_Y(\Delta_0)+c+d+e+\TIMETEMP(2(\bar\eta+c+d)+1)).
    \end{align}
    Then (\ref{Lemmastar2}) implies
    \begin{align}
        \frac{\P(\mc R_{k+1})}{\P(\mc R_k)} &\geq 1-\frac{\E[ Y_\gamma(\mathbf x_{k+1})|\mc R_k]}{\alpha_W(\Delta)}\\
        &> 1 - \frac{(1- p)^{-1}(\overline\alpha_Y(\Delta_0)+c+d+e+\TIMETEMP(2(\bar\eta+c+d)+1))}{\alpha_W(\Delta)}\\
        &\geq 1 - (1-(1- p)^{1/\TIMETEMP})\\
        &=(1- p)^{1/\TIMETEMP},
    \end{align}
    which concludes the induction step.\\
    \underline{Part 5: Conclusion.} It follows that
    \begin{align}\frac{\P(\mc R_\TIMETEMP)}{\P(\mc R_0)} = \frac{\P(\mc R_1)}{\P(\mc R_0)}\frac{\P(\mc R_2)}{\P(\mc R_1)}\dots\frac{\P(\mc R_\TIMETEMP)}{\P(\mc R_{\TIMETEMP-1})}>((1- p)^{1/\TIMETEMP})^\TIMETEMP = 1- p.\end{align}
    Since $\P(\mc R_0)=1$, we have $\P(\mc R_\TIMETEMP)>1- p$.
}

\section{Proof of Theorem \ref{thm:Recurrence}}\label{s:proof:thm:Recurrence}
\underline{Part 1: Preliminaries.} To prove Theorem \ref{thm:Recurrence}, based on Assumption \ref{assm:cost_controllability_detectability} we first define
\begin{flalign*}
     Y_\gamma(x)&:= V_\gamma(x)+\frac1\gamma W(x)-\frac{\gamma c_\gamma}{1-\gamma}+d+e_1&& \forall x\in\mc X,\\
    \overline\alpha_Y(s)&:=\overline\alpha_V(s)+2\overline\alpha_W(s)&&\forall s\in\R_{\geq0},\\
    e &:= e_1+e_2.
\end{flalign*}

The functions $ Y_\gamma$ serve a role similar to stochastic Lyapunov functions for recurrence considered in \cite{grammatico_discrete-time_2013}.
The following lemma provides lower and upper bounds for $ Y_\gamma$.

\begin{lem}\label{lem:Sandwich_bounds}
     Let Assumptions \ref{assum:measurable} and \ref{assm:cost_controllability_detectability} hold. Then, for any $\gamma\in[1/2,1)$ and $x\in\mc X$,
    \begin{align*}
        \alpha_W(\sigma(x)) \leq  Y_\gamma(x) \leq \overline\alpha_Y(\sigma(x)) + c + d + e.
    \end{align*}
\ediamond\end{lem}
\begin{pf}
    Let $\gamma\in[1/2,1)$ and $x\in\mc X$. To show the lower bound, let $\varepsilon>0$. By \eqref{eq:value_function}, there exists $h\in\mc H_a$ such that, for $u:=h(x)$,
    \begin{align*}
        V_\gamma(x)&\geq J_\gamma(x,h) - \varepsilon\\
        &= \ell(x,u) + \gamma\E J_\gamma(f(x,u,\mathbf v),h) - \varepsilon\\
        &\geq \ell(x,u) + \gamma\E V_\gamma(f(x,u,\mathbf v)) - \varepsilon.
    \end{align*}
    Using $W\geq0$ and $\gamma\leq1$, followed by this inequality, inequality \eqref{eq:dissipation} of part \ref{item:assm:detectability} of Assumption \ref{assm:cost_controllability_detectability} and the lower bound in part \ref{item:assm:cost_controllability} of Assumption \ref{assm:cost_controllability_detectability}, yields
    \begin{align*}
         Y_\gamma(x) &= V_\gamma(x) + \frac1\gamma W(x) -\frac{\gamma c_\gamma}{1-\gamma}+d+e_1\\
        &\geq V_\gamma(x) + W(x) -\frac{\gamma c_\gamma}{1-\gamma}+d+e_1\\
        &\geq\ell(x,u) + \gamma\underbrace{\E V_\gamma(f(x,u,\mathbf v))}_{\geq\frac {c_\gamma}{1-\gamma}-e_1} - \varepsilon + \underbrace{\E W(f(x,u,\mathbf v))}_{\geq0}\\
        &~~~- \ell(x,u) + \alpha_W(\sigma(x)) - d - \frac{\gamma c_\gamma}{1-\gamma} + d + e_1\\
        &\geq \frac{\gamma c_\gamma}{1-\gamma} - \gamma e_1 - \varepsilon + \alpha_W(\sigma(x)) - \frac{\gamma c_\gamma}{1-\gamma} + e_1\\
        &\geq\alpha_W(\sigma(x)) - \varepsilon.
    \end{align*}
    Because this holds for every $\varepsilon>0$, it follows that $ Y_\gamma(x) \geq\alpha_W(\sigma(x))$.
    Using the upper bounds in Assumption \ref{assm:cost_controllability_detectability}, as well as $c_\gamma\leq c$ and $\gamma\geq1/2$,
    \begin{align*}
         Y_\gamma(x) &= V_\gamma(x) + \frac1\gamma W(x) - \frac {c_\gamma}{1-\gamma}+c_\gamma+d+e_1\\
        &\leq \frac{c_\gamma}{1-\gamma}+e_2+\overline\alpha_V(\sigma(x)) + \frac1\gamma\overline\alpha_W(\sigma(x))\\
        &~~~-\frac{c_\gamma}{1-\gamma}+ c_\gamma+d+e_1\\
        &\leq \overline\alpha_V(\sigma(x)) + 2\overline\alpha_W(\sigma(x)) + c+d+e\\
        &= \overline\alpha_Y(\sigma(x)) + c+d+e.
    \end{align*}
    This completes the proof of the lemma.
    \mbox{}\hfill$\Box$
\end{pf}


\underline{Part 2: Setup.}
    Given $\Delta_0,\delta>0$ and $ p>0$, define
    \begin{align}\TIMETEMP:=\left\lceil\frac{4(\overline\alpha_Y(\Delta_0)+c+d+e)}{ p\delta}\right\rceil, \label{eq:proof:defJ}
    \end{align}
    where $\lceil\cdot\rceil$ denotes the ceiling function.
    By Lemma \ref{lem:Boundedness_with_high_prob}, there exists $\Delta>0$ such that $\Delta\geq\Delta_0$ and
    \begin{align}\forall\gamma\in(0,1)\,\forall x_0\in\sigma^{-1}(\Delta_0\B)\,\forall \mathbf x\in\mc S_\gamma^\eta(x_0): \P(\mc R_\TIMETEMP)>1-\frac p2,\label{eq:Thm_Recurrence_Delta}\end{align}
    where $\mc R_k$ for $k\in\Z_{\geq0}$ is defined as the event
    \begin{align*}
    		\mc R_k := \{&\omega\in\Omega~|~\forall j\in\{0,\dots,k\}:\sigma(\mathbf x_j(\omega))<\Delta\}.
    \end{align*}
    Note that $\mc R_k$ is measurable because $\mathbf x_0,\dots,\mathbf x_k$ and $\sigma$ are measurable.
    Now, choose $\gamma^\star\in(1/2,1)$ such that
    \begin{align} 
        (\gamma^\star)^{-1}(1-\gamma^\star)(\overline\alpha_V(\Delta)+e_2)\leq\delta/2.\label{eq:proof:gamma^*_def}
    \end{align}
    Let $\gamma\in(\gamma^\star,1)$, $x_0\in\sigma^{-1}(\Delta_0\B)$ and $\mathbf x\in\mc S_\gamma^\eta(x_0)$.
    In addition, let $\mathbf u_k, k\in\Z_{\geq0}$ be a corresponding sequence of controls as in Definition \ref{defn:solution}.
    For any $k\in\Z_{\geq0}$, define the event
    \begin{align*}
    \mc S_k := \{&\omega\in\Omega~|~\forall j\in\{0,\dots,k\}:\mathbf x_k(\omega)\notin\mc A_{\gamma,\delta}\}.
    \end{align*}
    Similar to before, $\mc S_k$ is measurable because $\mathbf x_0,\dots,\mathbf x_k, \sigma$ and $\mc A_{\gamma,\delta}$ are measurable.
    Note that
    \begin{align}
        \mc R_{k+1}\subseteq\mc R_k, \mc S_{k+1}\subseteq\mc S_k\quad\forall k\in\Z_{\geq0}. \label{eq:Rk_Sk_inclusions}
    \end{align}
    

    \noindent\underline{Part 3: Start of indirect proof.}
    We need to show $\P(\mc S_\TIMETEMP) < p$.
    Seeking a contradiction, assume $\P(\mc S_\TIMETEMP)\geq p$.
    Then, for any $k\in\{0,\dots,\TIMETEMP\}$, due to \eqref{eq:Thm_Recurrence_Delta} and \eqref{eq:Rk_Sk_inclusions},
    \begin{align}
        \P(\mc S_k\cap\mc R_k) &= \P(\mc S_k) + \P(\mc R_k) - \P(\mc S_k\cup\mc R_k)\\
        &\geq\P(\mc S_\TIMETEMP) + \P(\mc R_\TIMETEMP) - 1\\
        &> p + (1- p/2) - 1 =  p/2.\label{Sk_Rk_intersection_probability}
    \end{align}
    In particular, $\P(\mc S_k\cap\mc R_k)\neq0$.
    Thus, for any integrable random variable $\mathbf z:\Omega\to\R$ the conditional expectation $\E[\mathbf z|\mc S_k\cap\mc R_k]$, which we consider in subsequent parts, exists.
    The rest of this proof is dedicated to obtaining a contradiction under the assumption $\P(\mc S_T)\geq p$.

    \noindent\underline{Part 4: Derivation of bounds along the closed loop.}
    
        First, we show that $ Y_\gamma$, when conditioned on $\mc R_k\cap\mc S_k$, decreases by at least $\delta/2$ from $\mathbf x_k$ to $\mathbf x_{k+1}$. 
        For that, we use the strict causality guaranteed by Definition \ref{defn:solution}, as well as Assumption \ref{assm:cost_controllability_detectability}. 
    Due to strict causality, $\mathbf x_k$ and $\mathbf u_k$ are $\mc F_{k-1}$-measurable. Since $\mathbf v_j, j\in\Z_{\geq0}$ are i.i.d., $\mathbf v_k$ is independent of $\mc F_{k-1}$, and thus independent of $\mathbf x_k$ and $\mathbf u_k$.
    As $\mc R_k$ and $\mc S_k$ are $\mc F_{k-1}$-measurable as well, $\mathbf v_k$ is also independent of $\mc R_k\cap\mc S_k$.
    Thus, for any measurable $A\subseteq\mc X\times\mc U$ and $B\subseteq\mc V$,
    \begin{align*}
    &\P[(\mathbf x_k,\mathbf u_k,\mathbf v_k)\in A\times B|\mc R_k\cap\mc S_k]\\
    &=\P(\mc R_k\cap\mc S_k)^{-1}\P\left[((\mathbf x_k,\mathbf u_k)\in A) \wedge (\mathbf v_k\in B) \wedge (\mc R_k\cap\mc S_k)\right]\\
    &=\P(\mc R_k\cap\mc S_k)^{-1}\P\left[((\mathbf x_k,\mathbf u_k)\in A) \wedge (\mc R_k\cap\mc S_k)\right]\P[\mathbf v_k\in B]\\
    &=\P[(\mathbf x_k,\mathbf u_k)\in A|R_k\cap\mc S_k] \P[\mathbf v_k\in B].
    \end{align*}
    
    This means that the conditional push-forward distribution $\P^{\mathbf x_k,\mathbf u_k,\mathbf v_k|\mc R_k\cap\mc S_k}$ (that is, the probability measure on $\mc X\times\mc U\times\mc V$ given by
    $\P^{\mathbf x_k,\mathbf u_k,\mathbf v_k|\mc R_k\cap\mc S_k}(A) = \P[(\mathbf x_k,\mathbf u_k,\mathbf v_k)\in A|\mc R_k\cap\mc S_k]$
	for any measurable $A\subseteq\mc X\times\mc U\times\mc V$) can be decomposed as
    \begin{align*}
        \P^{\mathbf x_k,\mathbf u_k,\mathbf v_k|\mc R_k\cap\mc S_k} = \P^{\mathbf x_k,\mathbf u_k|\mc R_k\cap\mc S_k}\otimes\P^{\mathbf v_k},
    \end{align*}
    where $\otimes$ denotes the product measure.
    Using this along with Fubini's theorem, we can write any conditional expectation on $\mc R_k\cap\mc S_k$ as two nested integrals, the outer with respect to $\P^{\mathbf x_k,\mathbf u_k|\mc R_k\cap\mc S_k}(x,u)$, and the inner being the expectation with respect to $\mathbf v_k$ given fixed $(x,u)$.
    For brevity, we write $\P_k:=\P^{\mathbf x_k,\mathbf u_k|\mc R_k\cap\mc S_k}$ and $\E_k := \E[\cdot|\mc R_k\cap\mc S_k]$.
    This yields
    \begin{align*}
        &\E_k[ Y_\gamma(\mathbf x_{k+1})] - \E_k[ Y_\gamma(\mathbf x_k)]\\
        &= \int \E Y_\gamma(f(x,u,\mathbf v_k)) -  Y_\gamma(x)\,\td\P_k(x,u)\\
        &=\int \E V_\gamma(f(x,u,\mathbf v_k))-V_\gamma(x)\\
        &~~~~~~ +\gamma^{-1}(\E W(f(x,u,\mathbf v_k))-W(x))\,\td\P_k(x,u)
    \end{align*}
    Because $\mathbf u_k(\omega)\in\mc U_\gamma^\eta(\mathbf x_k(\omega))$ for all $\omega\in\Omega$, $V_\gamma(x) + \eta(x)\geq \ell(x,u) + \gamma\E V_\gamma(f(x,u,\mathbf v_k))$ for $\P_k$-almost every $(x,u)$ by \eqref{eq:U_gamma_eta}.
    Using this as well as part \ref{item:assm:detectability} of Assumption \ref{assm:cost_controllability_detectability} on the previously derived equality, the stage cost $\ell(x,u)$ cancels out and
    \begin{align*}
        &~~~~\E_k[ Y_\gamma(\mathbf x_{k+1})] - \E_k[ Y_\gamma(\mathbf x_k)]\\
        &\leq \int \gamma^{-1}((1-\gamma)\underbrace{V_\gamma(x)}_{\mathclap{\leq\overline\alpha_V(\Delta)+e_2+c_\gamma/(1-\gamma)}} +\, \eta(x) - \underbrace{\alpha_W(\sigma(x))}_{\mathrlap{\geq c_\gamma+d+\eta(x)+\delta}} +\, d)\,\td\P_k(x,u)\\
        &\leq\gamma^{-1}\left((1-\gamma)(\overline\alpha_V(\Delta)+e_2) - \delta\right)\\
        &\leq(\gamma^\star)^{-1}(1-\gamma^\star)(\overline\alpha_V(\Delta)+e_2) - \delta\stackrel{\eqref{eq:proof:gamma^*_def}}\leq -\delta/2,\tag{$\star_1$}\label{Thm_Recurrence_star1}
    \end{align*}
    where we used part \ref{item:assm:cost_controllability} of Assumption \ref{assm:cost_controllability_detectability} and the conditioning on $\mc R_k\cap\mc S_k$ from the second to third line, and $\gamma^\star\leq\gamma\leq1$ from the third to fourth line.

    We now move from $\mc R_k\cap\mc S_k$ to $\mc R_{k+1}\cap\mc S_{k+1}$.
    Since $ Y_\gamma(x)\geq0$ for every $x\in\mc X$ by Lemma \ref{lem:Sandwich_bounds}, it follows that for any $k\in\{0,\dots,\TIMETEMP-1\}$,
    \begin{align*}
        \E_{k+1}[ Y_\gamma(\mathbf x_{k+1})] &= \frac1{\P(\mc R_{k+1}\cap\mc S_{k+1})}\int_{\mc R_{k+1}\cap\mc S_{k+1}}\hspace{-0.2cm}  Y_\gamma(\mathbf x_{k+1})\td \P\\
        &\leq \frac1{\P(\mc R_{k+1}\cap\mc S_{k+1})}\int_{\mc R_k\cap\mc S_{k}}  Y_\gamma(\mathbf x_{k+1})\td \P\\
        &= \frac{\P(\mc R_k\cap\mc S_k)}{\P(\mc R_{k+1}\cap\mc S_{k+1})} \E_k[ Y_\gamma(\mathbf x_{k+1})] \tag{$\star_2$}\label{Thm_Recurrence_star2}
    \end{align*}
    (note that we are not dividing by zero thanks to \eqref{Sk_Rk_intersection_probability}).
    Because $\sigma(x_0)<\Delta_0\leq\Delta$, we have $\mc R_0=\Omega$ and, using the upper bound of Lemma \ref{lem:Sandwich_bounds},
    \begin{align}
        \E_0[ Y_\gamma(\mathbf x_0)] =  Y_\gamma(x_0) < \overline\alpha_V(\sigma(x_0))+c+d+e.\tag{$\star_3$}\label{Thm_Recurrence_star3}
    \end{align}
    The assumption $\P(\mc S_\TIMETEMP)\geq p$ of the indirect proof implies $\mc S_0 = \Omega$, because otherwise $x_0\in\mc A_{\gamma,\delta}$, implying $\P(\mc S_\TIMETEMP)=0$.

    \noindent\underline{Part 5: Conclusion.}
    Starting with \eqref{Thm_Recurrence_star3} and repeatedly applying \eqref{Thm_Recurrence_star1} and \eqref{Thm_Recurrence_star2} yields
    \begin{align*}
        &\E[ Y_\gamma(\mathbf x_\TIMETEMP)|\mc S_\TIMETEMP\cap R_\TIMETEMP]\\
        \begin{split}
        &< \underbrace{\frac{\P(\mc S_0\cap\mc R_0)}{\P(\mc S_\TIMETEMP\cap\mc R_\TIMETEMP)}}_{\mathclap{\leq 1/(p/2)\text{ by \eqref{Sk_Rk_intersection_probability}}}}(\overline\alpha_Y(\Delta_0)+c+d+e) \\
        &~~~~- \underbrace{\left(\frac{\P(\mc S_{0}\cap\mc R_{0})}{\P(\mc S_\TIMETEMP\cap\mc R_\TIMETEMP)} + \dots + \frac{\P(\mc S_{\TIMETEMP-1}\cap\mc R_{\TIMETEMP-1})}{\P(\mc S_\TIMETEMP\cap\mc R_\TIMETEMP)}\right)}_{\geq \TIMETEMP \text{ since } \mc S_0\cap\mc R_0\supseteq\dots\supseteq \mc S_\TIMETEMP\cap\mc R_\TIMETEMP}\frac\delta2 \stackrel{\eqref{eq:proof:defJ}}\leq 0,
        \end{split}
    \end{align*}
    which is a contradiction to $Y_\gamma(x)\geq0~\forall x\in\mc X$ by the lower bound of Lemma \ref{lem:Sandwich_bounds}.
    Thus, $\P(\mc S_\TIMETEMP)< p$, and
    \begin{align*}
        \P[\exists k\in\{0,\dots,\TIMETEMP\}:\mathbf x_k\in\mc A_{\gamma,\delta}] = 1-\P(\mc S_\TIMETEMP) > 1- p,
    \end{align*}
    concluding the proof.

\flushend

\end{document}